\documentclass[12pt]{article}
\usepackage{amssymb,amsmath}
\newtheorem{thm}{Theorem}
\newtheorem{cor}[thm]{Corollary}
\newtheorem{lemma}[thm]{Lemma}
\newenvironment{defin}{\medskip\noindent{\sc
Definition}.}{\goodbreak\medskip}

\newenvironment{remk}{\noindent{\sc
Remark}.}{\goodbreak\vskip10pt}

\def\demo{\medskip\goodbreak\noindent
     \hbox{\sc Proof \kern .3em}\ignorespaces}%
  \def \qedbox{$\square$}%
  \def \qed{\hglue1mm\hfill{\ifmmode\qedbox
     \else\unskip\ \hglue0mm\hfill\qedbox\medskip
      \goodbreak\fi}}%
\def\enddemo{\qed\goodbreak\vskip10pt}%

\def\qed{\hglue1mm\hfill\raise -2pt\hbox{\vrule\vbox to 10pt{\hrule width
4pt
                  \vfill\hrule}\vrule}}

\newcommand{\T}{\mathbb {T}}

\newcommand{\R}{\mathbb {R}}

\newcommand{\N}{\mathbb {N}}

\newcommand{\Vc}{\mathcal {V}}

\newcommand{\Ec}{\mathcal {E}}

\newcommand{\Gc}{\mathcal {G}}

\newcommand{\Ac}{\mathcal {A}}

\begin{document}
\title{{
}
Pseudographs and Lax-Oleinik semi-group~: a geometric and dynamical  interpretation }
\author{M.-C. ARNAUD
\thanks{ANR KAM faible}
\thanks{Universit\'e d'Avignon et des Pays de Vaucluse, Laboratoire d'Analyse non lin\' eaire et G\' eom\' etrie (EA 2151),  F-84 018Avignon,
France. e-mail: Marie-Claude.Arnaud@univ-avignon.fr}
}
\maketitle
\abstract{Let $H~: T^*M\rightarrow \R$ be a Tonelli Hamiltonian defined on the cotangent bundle of a compact and connected manifold and let $u~: M\rightarrow \R$ be a semi-concave function. If $\Ec(u)$ is the set of all the super-differentials of $u$ and $(\varphi_t)$ the Hamiltonian flow of $H$, we prove that for $t>0$ small enough, $\varphi_{-t}(\Ec (u))$ is an exact Lagrangian Lipschitz graph.\\ 
This provides a geometric interpretation/explanation of a regularization tool that was introduced by  P.~Bernard in \cite{Be2} to prove the existence of $C^{1,1}$ subsolutions.

}
\newpage
\tableofcontents
\newpage
 
\section{Introduction} 
In the recent developments  of the so-called ``weak K.A.M. theory'', the notion of ``pseudograph'' did appear recently in an article of P.~Bernard (see \cite{Be1}) to prove some results concerning Arnold's and Mather's diffusion.  Let us explain quickly  how this notion appeared. \\
We consider the Hamilton-Jacobi equation~: $H(x, du(x))=C$ for a Hamiltonian function $H~: T^*M\rightarrow \R$ defined on a cotangent bundle that is $C^2$, superlinear and convex in the fiber. In the 1980's, P.-L.~Lions and M.~Crandall introduced the notion of viscosity solution for this equation (see \cite{Cra-Li}). In the case   $M=\T^n$, Lions, Papanicolaou and Varadhan proved  the existence of a viscosity solution. Then, in \cite{Fa1}, A.~Fathi proved the existence of  a viscosity solution (that he called a weak K.A.M. solution) for any manifold.  Such a weak K.A.M. solution is semi-concave and hence locally Lipschitz (see for example \cite{Fa}). A semi-concave function $u~: M\rightarrow \R$ is   Lipschitz and hence differentiable on a set $E\subset M$ with full (Lebesgue) measure, and the graph $\Gc (u)=\{ (q, du(q)); q\in E\}$ of the derivatives of any semi-concave function is what we call  a {\em pseudograph}. \\

When $u$ is $C^2$, the pseudograph $\Gc (u)$ is in fact a graph above the whole manifold $M$ and is a Lagrangian graph. 
\medskip

That's why a very natural question is~: 

\medskip
\noindent{\bf Questions ~:} are the pseudographs Lagrangian manifolds in general? And, as a pseudograph is not a smooth manifold,  in which sense?\\
\medskip

Let us notice that in the other sense, M.~Chaperon proved in \cite{Cha} that every  Lagrangian submanifold of $T^*M$ that is Hamiltonianly isotopic to the zero section can be ``cut'' in such a way that we obtain the graph of the differential of a Lipschitz function defined on $M$. In some cases, A.~ÊOttolenghi \& C.~Viterbo proved in \cite{Ot-Vi} that this Lipschitz function is a semi-concave one,  and hence the ``cut graph'' is  a pseudograph. Let us notice too that we proved in \cite{Arna2} that   any invariant Lagrangian manifold that is Hamiltonlianly isotopic to the zero section and invariant by a Tonelli Hamiltonian is the graph of a smooth function. Hence if the pseudograph of a weak KAM solution is obtained by cutting an invariant Lagrangian submanifold that is Hamiltonianly isotopic to the zero section, then this pseudograph has to be a true smooth submanifold.

\medskip
Once we have proved that the pseudographs are some ``Lagrangian manifolds'' (in some sense that we will explained soon), we know that their images by the Hamiltonian flows are Lagrangian too because a Hamiltonian flow is symplectic. But in general the image of a pseudograph by a Hamiltonian flow is not a pseudograph (it may happen that it is not a graph).  To stay in the class of the graphs, let us consider  the two Lax-Oleinik semi-groups $T_t, \breve T_t~: C^0(M, \R)\rightarrow C^0(M, \R)$ associated to the considered Tonelli Hamiltonian (they will be precisely defined). Let us recall some well-known results concerning the relationships between the action of these two semi-groups and the action of the Hamiltonian flow on the pseudographs (see \cite{Fa} and \cite{Be1}). We denote the associated Hamiltonian flow by $(\varphi_t)$.
\begin{enumerate}
\item for every $t>0$, all the functions of $T_t(C^0(M, \R))$ (resp. $\breve T_t(C^0(M, \R))$) are semi-concave (resp. semi-convex);
\item if $u~: M\rightarrow \R$ is a semi-concave function, then for all $t>0$, we have $\overline{\Gc(T^tu)}\subset \varphi_t(\Gc(u))$; the action of the negative Lax-Oleinik semi-group on the derivative of a semi-concave function is what follows~: we take the image $\varphi_t(\Gc (u))$ of the graph of $du$ by the positive flow and we remove some part of this set;
\item similarly, if $u$ is semi-convex, we have for all $t>0$~: $\overline{\Gc(\breve T^tu)}\subset \varphi_{-t}(\Gc(u))$;
\item if we just assume that $u$ is continuous, then for all $t>0$ the set $\varphi_{-t}\left ( \overline{\Gc(T^tu)}\right)$ is a subset of the set of the sub-derivatives of $u$ and $\varphi_t\left( \overline{\breve \Gc(T^tu)}\right)$ is a subset of the set of the super-derivative of $u$. Hence the positive Lax-Oleinik semi-group maps any continuous function on a function $\breve T_tu$ such that $\Gc(\breve T_tu)$ is  a part of the image by the negative flow of what we will call the enlarged pseudograph (i.e. the set of all the super-derivatives of $u$).
\end{enumerate}
Hence there is a deep link between the action of the Lax-Oleinik semi-group on the graphs of  super/sub-derivatives and the action of the Hamiltonian flow. Our purpose is to give a precise statement concerning the action of the positive  Lax-Oleinik semi-group $(\breve T_t)$ on the semi-concave functions and to prove simultaneously that the enlarged pseudographs of the semi-concave functions are some Lipschitz Lagrangian submanifolds.

\medskip

Before explaining our result,   let us introduce precisely some notions. At first, we recall what is a semi-concave function and we define the enlarged pseudographs.

\begin{defin} \begin{enumerate}
\item Let $U$ be an open subset of $\R^d$, $K\geq 0$ be a constant and $u~: U\rightarrow \R$ be a function. We say that $u$ is $K$-semi-concave if for every $x\in U$, there exists a linear form $p_x$ defined on $\R^d$ such that ~:
$$\forall y\in U, u(y)\leq u(x)+p_x(y-x)+K\| y-x\|^2;$$
(where $\| .\|$ is the usual Eucidian norm). Then we say that $p_x$ is a $K$-super-differential of $u$ at $x$.
\item Let $M$ be a compact and connected manifold with a finite atlas $\Ac=\{ (U_i, \Phi_i~: U_i\rightarrow  \R^d); 1\leq i\leq N\}$ and $u~: M\rightarrow \R$ be a function; we say that $u$ is $K$-semi-concave if for every $i\in \{ 1, \dots , N\}$, the function $u\circ \Phi_i^{-1}~: \Phi_i(U_i)\rightarrow \R$ is $K$-semi-concave. Then, a $K$-superdifferential of $u$ is a $p_x\circ D\Phi_i(x)$ where $p_x$ is a $K$-super-differential of $u\circ \Phi_i^{-1}$ at $\Phi_i(x)$.
\item A function is semi-concave if it is   $K$-semi-concave for a certain $K$; while the quantitative notion of ``$K$-semi-concave function'' depends on the considered atlas of $M$ that we choose, the notion of ``semi-concave function'' is independent of this atlas. The notion of super-differential too doesn't depend on the atlas.
\item if $u~: M\rightarrow \R$ is semi-concave, its {\em enlarged pseudograph} is the set $\Ec (u)$ of all the super-differentials of $u$~:
$$\Ec(u)=\{ (x, p_x); p_x\quad{\rm is}\quad{\rm a} \quad{\rm superdifferential}\quad{\rm of}\quad  u \quad{\rm at}\quad x\}.$$
\end{enumerate}
\end{defin}
The enlarged pseudograph $\Ec (u)$ of a semi-concave function $u$ contains its pseudograph $\Gc (u)$;  in general, $\Ec(u)$ is no longer a graph and $\Ec(u)$ is compact (it's clearly closed and P.~Bernard proved in \cite{Be1} that it is bounded).\\

\begin{remk}
In fact, even if it doesn't appear in the notation,  the definition of $\Ec (u)$  depends on the choice of the constant $K$ of semi-concavity that we choose, and in the proofs we will fix such a constant $K$. But a posteriori, because of theorem \ref{T1}, we see that $\Ec (u)$ is independant of this constant.
\end{remk}

For a survey of the principal properties of the semi-concave functions, the reader may have a look at the appendix of \cite{Be1} and the book \cite{Fa}.\\

Let us now explain which kind of submanifolds will interest us~:

\begin{defin} Let $M$ be a $d$-dimensional compact and connected manifold.
 \begin{enumerate}
\item a non-empty subset $N$ of $T^*M$ is a $d$-dimensional {\em Lipschitz submanifold} of $T^*M$   if  for   every $x\in N$, there exists a (smooth) chart $(U, \Phi)$ of $T^*M$ at $x$ such that $\Phi (N\cap U)$ is the graph of a Lipschitz map $\ell~: V\rightarrow \R^d$ defined on a open subset $V$ of $\R^d$. Of course, this notion is invariant by  $C^1$-diffeomorphism.
\item a Lipschitz graph is $\{ s(x); x\in M\}$ where $s~: M\rightarrow T^*M$ a Lipschitz section. Of course, a Lipschitz graph is a $d$-dimensional Lipschitz submanifold of $T^*M$. 
\item a $d$-dimensional Lipschitz submanifold $N$ of $T^*M$ is {\em exact Lagrangian} if it is exact Lagrangian in the sense of distributions, that is if for every $\gamma~: [a, b]\rightarrow  N$ closed Lipschitz arc drawn on $N$, we have~: $0= \int_\gamma\lambda$ (where $\lambda$ designates the Liouville $1$-form of $T^*M$). This notion is invariant under $C^1$ exact symplectic diffeomorphisms. 
\end{enumerate}
\end{defin} 
Then the  Lipschitz graph of a  Lipschitz section $s~: M\rightarrow T^*M$ is exact Lagrangian if and only if there exists a $C^{1, 1}$ function $u~: M\rightarrow \R$ (that is a $C^1$ function whose derivative is Lipschitz) such that~: $s=du$.

The result that we obtain is~:
\begin{thm}\label{T1}
Let $M$ be a  compact  and connected manifold,  let $u~: M\rightarrow \R$ be a semi-concave function and let $\Ec (u)$ be its enlarged pseudograph. Let $(\varphi_t)$ be a Tonelli Hamiltonian flow of $T^*M$.Then there exists $\varepsilon>0$ such that for all $t\in ]0, \varepsilon]$, we have~: $\Gc (d\breve T_tu)=\varphi_{-t}(\Ec (u))$ is a Lipschitz graph above the whole manifold.
\end{thm}
We immediately deduce~:
\begin{cor}
The enlarged pseudograph of any semi-concave function of $M$ is a Lipschitz exact Lagrangian submanifold of $T^*M$. 
\end{cor}
Let us recall that in \cite{Be2}, P.~Bernard proved the following regularization result (that he used to prove the existence of $C^{1,1}$ sub-solutions)~:  for each semi-concave function $u~: M\rightarrow \R$, for every $t>0$ small enough, the function $\breve T_tu$ is $C^{1,1}$. Of course we reprove this result, but this is not the goal of this article and Bernard's proof is shorter and more efficient. Our purpose is to give a geometric/dynamical interpretation in terms of exact Lagrangian Lipschitz sub-manifold and in term of Hamiltonian flow of the action of the Lax-Oleinik semi-group on the enlarged pseudo-graphs.

\section{Proof of theorem \ref{T1}}

We assume that  $M$ is a compact and connected manifold with a finite atlas $\Ac$ and $u~: M\rightarrow \R$ is a $K$-semi-concave function. We consider any Tonelli Hamiltonian function $H~: T^*M\rightarrow \R$ and denote by $(\varphi_t)_{t\in\R}$ its Hamiltonian flow. 
\subsection{Proof that $\varphi_t(\Ec(u))$ is a graph for $t\in [-\varepsilon, 0[$}\label{SS21}

Given $\varepsilon\in ]0, 1]$ small enough, we want to know if it is possible that for a $t\in [-\varepsilon, 0[$, $\varphi_t(\Ec (u))$ is not a graph above  a certain part of $M$. 
\\
To prove that, we will need some inequalities given in the following lemmata.

  \begin{lemma}\label{L3}
  We assume that $(q_0, p_0), (q_1, p_1)\in \Ec(u)$ are in a same chart of the atlas. Then~:
   $$(p_1-p_0)(q_1-q_0)\leq 2K\| q_1-q_0\|^2.$$
  \end{lemma}
  
 \demo We know too that $p_j$ is a $K$-super-differential  of $u$ at $q_j$, for $j=0,1$. Hence~:\\
  $u(q_1)-u(q_0)\leq p_0(q_1-q_0)+K\| q_1-q_0\|^2$;\\
  $u(q_0)-u(q_1)\leq p_1(q_0-q_1)+K\| q_1-q_0\|^2$.\\
  By adding up these two inequalities, we deduce the lemma.\enddemo
  \begin{lemma}\label{L4} Let $K$ be a compact subset of $T^*M$ that is convex in the fiber  and let $c, C$ be two constants such that~: $$\forall \tau \in [-1, 1], \forall x\in K, \forall v\in\R^d, c\| v\|^2\leq H_{p, p}(\varphi_\tau (x))(v,v)\leq C\| v\|^2.$$ Then there exists $\varepsilon>0$ such that, for every $t\in ]0, \varepsilon]$  and every $(q, p)$, $(q, p+\Delta p)\in K$, if we use the notations~: $(q_0, p_0)=\varphi_t(q, p)$ and $(q_1, p_1)=\varphi_t(q, p+\Delta p)$, we have~:
  $$(p_1 -p_0 )(q_1 -q_0 )\geq \frac{c}{2}t\| \Delta p\|^2\quad {\rm  and}\quad \|q_1-q_0\|\leq2Ct\| \Delta p\|.$$
  
 \end{lemma}
\demo
Because $K$ is compact, if we choose $\varepsilon>0$ small enough, then $q_0$, $q_1$ and $\pi\circ \varphi_\tau (q_0, p_0)_{\tau\in [-\varepsilon, 0]}$ are in a same chart of the atlas $\Ac$.\\
 Then from now we work in the coordinates given by such a chart, i.e. in $\R^d$ and $T^*\R^d=\R^d\times \R^d$ and we write~: $(q_1, p_1)\-=\varphi_t(q, p+\Delta p)$ with $t\in ]0, \varepsilon]$.\\
 As $K$ is compact and $t\in [0, 1]$, there exists a constant $R>0$ such that, necessarily~: $\| \Delta p\| \leq R$ (for the usual Euclidian norm in $\R^d$). \\
 We compute (let us notice that for every $s\in [0, 1]$, we have~: $(q , p +s\Delta p)\in K$ because $K$ is convex in the fibers)~:
$$\varphi_t(q, p+\Delta p)-\varphi_t(q, p)=\int_0^1D\varphi_t(q, p+s\Delta p)(0, \Delta p)ds.$$
Using the linearized Hamilton equations, we obtain that the quantity in the integral is equal to~:
 $$ (t(H_{p,p}(\varphi_t(q, p +s\Delta p)\Delta p+ \| \Delta p\| \varepsilon_1 (s, t, \Delta p)), \Delta p+  \|\Delta p\|\varepsilon_2(s,t, \Delta p))$$
 where the functions $\epsilon_j$ tend  uniformly to $0$ when $t$ tends to $0$.
 We deduce~:
 $$(p_1 -p_0 )(q_1 -q_0 )=t\int_0^1\left[ H_{p,p}(\varphi_t((q, p+s\Delta p)(\Delta p, \Delta p) +t\| \Delta p\|^2\eta (s,t, \Delta p)\right]ds$$
  where the function  $\eta $ tends  uniformly to $0$ when $t$ tends to $0$.
Hence if $\varepsilon$ has been chosen small enough, we have~: 
  $$(p_1 -p_0 )(q_1 -q_0 )\geq \frac{c}{2}t\| \Delta p\|^2\quad {\rm  and}\quad \|q_1-q_0\|\leq2Ct\| \Delta p\|.$$
 \\

\enddemo

If $x= (q, p)\in T^*M$, we denote by $\Vc (x)$ its vertical~: $\Vc (x)=T_q^*M=\{ y\in T^*M; \pi (y)=\pi (x)=q\}$ where $\pi~: T^*M\rightarrow M$ designates the usual projection. \\
Then we want to know  if it is possible for a $t\in ]0, \varepsilon]$ and a $x\in \Ec (u)$  that $\Vc(\varphi_{-t}(x))\cap \varphi_{-t}(\Ec(u))$ contains at least two points. It means that there exists two different points $(q_0, p_0), (q_1, p_1)\in \Ec(u)$ such that $(q_1, p_1)\in \varphi_t(\Vc (\varphi_{-t}(q_0, p_0)))$. We use the notation~: $(q, p)=\varphi_{-t}(q_0, p_0)$ and $(q, p+\Delta p)=\varphi_{-t}(q_1, p_1)$.\\
As $\Ec (u)$ is compact subset of $T^*M$ that is compact in the fiber, we can use lemma \ref{L3} to choose $\varepsilon>0$. Then we have~:
  $$(p_1 -p_0 )(q_1 -q_0 )\geq \frac{c}{2}t\| \Delta p\|^2\quad {\rm  and}\quad \|q_1-q_0\|\leq2Ct\| \Delta p\|.$$
  Now lemma \ref{L3} tells us that~: $(p_1-p_0)(q_1-q_0)\leq 2K\| q_1-q_0\|^2$.
  Then~:
  $$(p_1-p_0)(q_1-q_0)\leq 2K\| q_1-q_0\|^2\leq 8C^2Kt^2\|\Delta p\|^2.$$
  Finally, we have proved that there exist two strictly positive constants $c$ and $C$  such that~: 
   $$\frac{c}{2}t\| \Delta p\|^2\leq (p_1-p_0)(q_1-q_0) \leq 8C^2Kt^2\|\Delta p\|^2.$$
   It is obviously impossible for $t>0$ small enough and $\Delta p\not= 0$.
   
   \subsection{Proof that $\pi\circ \varphi_t(\Ec(u))=M$ }\label{SSGR}

 We want to prove that for $t\in [-\varepsilon, 0[$, the graph $\varphi_t(\Ec(u))$ covers the whole $M$.\\
 We have recalled in introduction that for all $t>0$, we have~: $\overline{\Gc (\breve T_tu)}\subset\varphi_{-t}(\Ec (u))$;  this implies directly that $\pi\circ \varphi_t(\Ec(u))=M$.

 \subsection{Proof that $\ \varphi_t(\Ec (u))$ is a Lipschitz graph}
 We have proved that for $t\in ]0, \varepsilon_0]$, $\ \varphi_{-t}(\Ec (u))$ is a graph above $M$. Because this graph is compact ($\Ec (u)$ is compact), it's the graph of a continuous section $s_t~: M\rightarrow T^*M$.\\
 We have to prove that $s_t$ is Lipschitz. We may eventually change $\varepsilon_0$ in such a way that $K<\frac{1}{C\varepsilon_0}$.\\
 We will use the so-called Bouligand's paratingent cone~: 
 
 \begin{defin}
 Let $E$ be a subset of $T^*M$. The paratingent cone to $E$ at $(q,p)\in E$ is defined (in chart but it doesn't depend on the chart) as the subset of $T_{(q,p)}(T^*M)$ whose elements are the limits of the sequences~:\\
  $\left( \frac{1}{t_n}(q_n-q'_n), \frac{1}{t_n}(p_n-p'_n)\right)_{n\in\N}$ with $t_n\in \R^*_+$, $q_n, q'_n, p_n, p'_n\in E$ and $\lim q_n=\lim q'_n=q$, $\lim p_n=\lim p'_n=p$. It is denoted by $C_{(q,p)}E$.
 \end{defin}
 If $(q,p), (q',p')\in\Ec (u)$ are in a same chart, we have proved in lemma \ref{L3}  that~:
 $(p'-p)(q'-q)\leq 2K\| q'-q\|^2$. We deduce that for all $(\delta q, \delta p)\in C_{(q,p)}\Ec (u)$, we have~: $\delta p.\delta q\leq 2K\| \delta q\|^2$.\\
 
  Moreover, we deduce easily from lemma \ref{L4} that if $R>0$, there exists $\varepsilon>0$ such that for every $(q, p), (q, p+\Delta p)\in T^*M$ that satisfy  $\| p\|\leq R$ and $\| p+\Delta p\| \leq R$, we have if we use the notations $\varphi_t(q,p)=(q_0, p_0)$ and $\varphi_t(q, p+\Delta p)=(q_1, p_1)$ for a $t\in ]0, \varepsilon]$~: 
 $$(p_1 -p_0 )(q_1 -q_0 )\geq \frac{c}{2}t\| \Delta p\|^2\quad {\rm  and}\quad \|q_1-q_0\|\leq2Ct\| \Delta p\|.$$
 Looking at what happens when $\Delta p$ tends to $0$, we deduce that~:\\
for every $(q, p)\in \varphi_{-t}(\Ec (u))$, for every $\delta p_0\in T_{(q,p)}(T_q^*M)$, if we use the notation $D\varphi_t(q,p)(0, \delta p_0)=(\delta q, \delta p)$, then we have~:  
 $$\delta p.\delta q\geq \frac{c}{2}t\| \delta p_0\|^2\quad {\rm  and}\quad \|\delta q\|\leq2Ct\| \delta p_0\|;$$
 then~: $  \| \delta q\|^2\leq 8\frac{C^2}{c}t \delta p.\delta q$. \\
 
 Finally, we have proved for $(q, p)\in \Ec(u)$ that~: 
 \begin{enumerate}
 \item[$\bullet$]   for all $(\delta q, \delta p)\in C_{(q,p)}\Ec (u)$, we have~: $\delta p.\delta q\leq 2K\| \delta q\|^2$;
 \item[$\bullet$] for all $(\delta q, \delta p)\in T_{(q, p)}^*M$ that is in the image by $D\varphi_t $ of the vertical $V(\varphi_{-t}(q,p))=\ker D\pi(\varphi_{-t}(q,p))$, we have~: $  \| \delta q\|^2\leq 8\frac{C^2}{c}t \delta p.\delta q$.
 \end{enumerate}
 If we choose $\varepsilon < \frac{Kc}{4C^2}$, we obtain that $D\varphi_tV(\varphi_{-t}(q,p))\cap C_{(q,p)}\Ec (u)=\{ 0\}$, and then that~:  $V(\varphi_{-t}(q,p))\cap  D \varphi_{-t}(C_{(q,p)}(\Ec (u)))=\{ 0\}$.

Finally, we have proved that the    paratingent cone to $\varphi_{-t}(\Ec(u))$, which is the graph of $s_t$,  contains no vertical line. Let us deduce that $s_t$ is Lipschitz. We assume that there are two sequences of points $(q_n, p_n)$, $(q'_n, p'_n)$  of $\varphi_{-t}(\Ec(u))$ such that $\lim \frac{\| p'_n-p_n\|}{d(q'_n,q_n)}=+\infty$.  Using a subsequence, because $\varphi_{-t}(\Ec (u))$ is compact, we may assume that the two sequences converge. Then necessarily $(q_n)$ and $(q'_n)$ have the same limit (because the previous limit is $+\infty$ and $\| p'_n-p_n\|$ is bounded). Hence by continuity of $s_t$,  $(p_n)$ and $(p'_n)$ too have the same limit. But if we write $t_n= \| p_n-p'_n\| $ and if we use a subsequence in such a way that $(\frac{p'_n-p_n}{\| p'_n-p_n\|})$ converges to a $u$, we obtain that $\displaystyle{\lim_{n\rightarrow \infty}\frac{1}{ t_n}(q'_n-q_n, p'_n-p_n)=(0, u)}$ is in the paratingent cone to $\varphi_{-t}(\Ec (u))$ at $(q,p)$, it  contradicts the fact that this paratingent cone contains no vertical line.  Hence $s_t$ is Lipschitz.

\subsection{Proof that $\ \varphi_t(\Ec (u))$ is an exact Lagrangian Lipschitz graph}
We have to prove that there exists a $C^1$ function (hence it will be $C^{1, 1}$) $u_t~: M\rightarrow \R$ such that $s_t=du_t$. It is enough to prove that for any closed Lipschitz arc $\gamma~: [a, b]\rightarrow M$, then $\int_a^bs_t(\gamma (\tau ))\dot\gamma(\tau)d\tau=0$. Let us define a  closed loop of $T^*M $ by~: $\forall \tau\in [a,b], \eta (\tau)=(\eta_1(\tau), \eta_2(\tau))=\varphi_t(\gamma(\tau), s_t(\gamma (\tau)))$. The arc $\gamma$ being Lipschitz and $s_t$ being Lipschitz, the arc $\eta$ is Lipschitz too. Hence we can define $\int_\eta\lambda $ where $\lambda$ is the Liouville 1-form. The flow being exact symplectic, we have~: $\int_\eta\lambda  =\int_a^bs_t(\gamma (\tau ))\dot\gamma(\tau)d\tau$. We are reduced to compute $\int_\eta\lambda=\int_a^b\eta_2(\tau)\dot\eta_1(\tau )d\tau$.\\
Let us recall that $\eta$ is a closed Lipschitz arc drawn on $\Ec (u)$; then $\eta_2(\tau)$ is a $K$-super-differential of $u$ at $\eta_1(\tau)$ and~:
$$u(\eta_1 (\tau+\delta \tau))-u(\eta_1(\tau))\leq \eta_2(\tau)(\eta_1(\tau+\delta \tau)-\eta_1(\tau)) +K\| \eta_1(\tau +\delta \tau)-\eta_1(\tau)\|^2.$$
Moreover, $u$, $\eta_1$ and $\eta_2$ are Lipschitz, then (Lebesgue) almost everywhere derivable. If $\tau$ is a point where $u\circ\eta_1$ and $\eta_1$ are derivable, we obtain by dividing by $\delta t$ (positive or negative) and taking the limit when $\delta t$ tends to $0$~:

$\frac{d{}}{dt}(u\circ\eta_1)(\tau)=\eta_2(\tau)\dot\eta_{1}(\tau)$ and by integration~:
$$\int_a^b\eta_2(\tau)\dot\eta_{1}(\tau)d\tau=\int_a^b\frac{d{}}{dt}(u\circ\eta_1)(\tau)d\tau=u(\eta_1(b))-u(\eta_1(a))=0.
$$
\subsection{Proof that $\Gc (d\breve T_tu)=\varphi_{-t}(\Ec (u))$}
We have proved that $\varphi_{-t}(\Ec (u))$ is an exact Lagrangian Lipschitz graph for every $t\in [-\varepsilon , 0[$; we write~: $\varphi_{-t}(\Ec (u))=\Gc(u_t)$ with $u_t$ that is $C^{1,1}$. Moreover, we have noticed in the introduction that $\Gc(\breve T_t u)\subset \varphi_{-t}(\Ec (u))=\Gc(u_t)$; as $\breve T_tu$ is Lipschitz, we deduce that for Lebesgue almost every $q\in M$, we have~: $du_t=d\breve T_tu$. The derivative of the two Lipschitz functions $u_t$ and $\breve T_tu$ are almost everywhere equal, then $\breve T_tu-u_t$ is a constant function and then $\breve T_t u$ is $C^{1,1}$ and we have the equality~: $\Gc (d\breve T_tu)=\varphi_{-t}(\Ec (u))$.

\end{document}